\documentclass[reqno]{amsart}
\usepackage{amssymb}
\usepackage[hypertex]{hyperref}
\allowdisplaybreaks[4]
\theoremstyle{plain}
\newtheorem{thm}{Theorem}
\theoremstyle{remark}
\newtheorem{rem}{Remark}
\DeclareMathOperator{\td}{d}
\date{Completed on 30 January 2010 in Tianjin}
\date{}

\begin{document}

\title{Sharp bounds for harmonic numbers}

\author[F. Qi]{Feng Qi}
\address[F. Qi]{Department of Mathematics, College of Science, Tianjin Polytechnic University, Tianjin City, 300160, China}
\email{\href{mailto: F. Qi <qifeng618@gmail.com>}{qifeng618@gmail.com}, \href{mailto: F. Qi <qifeng618@hotmail.com>}{qifeng618@hotmail.com}, \href{mailto: F. Qi <qifeng618@qq.com>}{qifeng618@qq.com}}
\urladdr{\url{http://qifeng618.spaces.live.com}}

\author[B.-N. Guo]{Bai-Ni Guo}
\address[B.-N. Guo]{School of Mathematics and Informatics, Henan Polytechnic University, Jiaozuo City, Henan Province, 454010, China}
\email{\href{mailto: B.-N. Guo <bai.ni.guo@gmail.com>}{bai.ni.guo@gmail.com},
\href{mailto: B.-N. Guo <bai.ni.guo@hotmail.com>}{bai.ni.guo@hotmail.com}}

\begin{abstract}
In the paper, we first survey some results on inequalities for bounding harmonic numbers or Euler-Mascheroni constant, and then we establish a new sharp double inequality for bounding harmonic numbers as follows: For $n\in\mathbb{N}$, the double inequality
\begin{equation*}
-\frac{1}{12n^2+{2(7-12\gamma)}/{(2\gamma-1)}}\le H(n)-\ln
n-\frac1{2n}-\gamma<-\frac{1}{12n^2+6/5}
\end{equation*}
is valid, with equality in the left-hand side only when $n=1$, where the scalars $\frac{2(7-12\gamma)}{2\gamma-1}$ and $\frac65$ are the best possible.
\end{abstract}



\thanks{The authors were supported in part by the Science Foundation of Tianjin Polytechnic University}

\thanks{This paper was typeset using \AmS-\LaTeX}

\maketitle

\section{Introduction}

The series
\begin{equation}
1+\frac12+\frac13+\dotsm+\frac1n+\dotsm
\end{equation}
is called harmonic series. The $n$-th harmonic number $H(n)$ for $n\in\mathbb{N}$, the sum of the first $n$ terms of the harmonic series, may be given analytically by
\begin{equation}\label{hn-anal}
H(n)=\sum_{i=1}^n\frac1i=\gamma+\psi(n+1),
\end{equation}
see \cite[p.~258, 6.3.2]{abram}, where $\gamma=0.57721566\dotsm$ is Euler-Mascheroni constant and $\psi(x)$ denotes the psi function, the logarithmic derivative $\frac{\Gamma'(x)}{\Gamma(x)}$ of the classical Euler gamma function $\Gamma(x)$ which may be defined by
\begin{equation}\label{gamma-dfn}
\Gamma(x)=\int^\infty_0t^{x-1} e^{-t}\td t,\quad x>0.
\end{equation}
\par
In \cite{PS}, the so-called Franel's inequality in literature was given by
\begin{equation}\label{franel}
\frac1{2n}-\frac1{8n^2}<H(n)-\ln n-\gamma<\frac1{2n},\quad n\in\mathbb{N}.
\end{equation}
In \cite[pp.~105--106]{klam}, by considering
\begin{equation}
I_n=\int_{1/n}^1\biggl(\frac1x-\biggl[\frac1x\biggr]\biggr)\td x=\ln n-H(n)
\end{equation}
and $0<I_n<\frac12$, where $[t]$ denotes the largest integer less tan or equal
to $t$, it was established that
\begin{equation}
\frac12<H(n)-\ln n<1,\quad n\in\mathbb{N}.
\end{equation}
In \cite[pp.~128--129, Problem~65]{klam}, it was verified that
\begin{equation}
\label{p65} \frac12\ln(2n+1)<\sum_{k=1}^n\frac1{2k-1}<1+\frac12\ln(2n-1),\quad n\in\mathbb{N}.
\end{equation}
\par
In \cite{young}, it was obtained that
\begin{equation}
\frac1{2(n+1)}<H(n)-\ln n-\gamma<\frac1{2n},\quad n\in\mathbb{N}.
\end{equation}
In \cite{DeTemple}, it was proved that
\begin{equation}
\frac1{24(n+1)^2}<H(n)-\ln\biggl(n+\frac12\biggr)-\gamma<\frac1{24n^2},\quad n\in\mathbb{N}.
\end{equation}
\par
In \cite{toth-3432}, the following problems were proposed:
\begin{enumerate}
\item
Prove that for every positive integer $n$ we have
\begin{equation}\label{toth-ineq}
\frac1{2n+2/5}<H(n)-\ln n-\gamma<\frac1{2n+1/3}.
\end{equation}
\item
Show that $\frac25$ can be replaced by a slightly smaller number, but
$\frac13$ cannot be replaced by a slightly larger number.
\end{enumerate}
In \cite{high-toth}, by using
\begin{equation}\label{hn-4}
H(n)=\ln n+\gamma+\frac1{2n}-\frac1{12n^2} +\frac{\varepsilon_n}{120n^4}
\end{equation}
for $0<\varepsilon_n<1$, these problems were answered affirmatively. The
editorial comment in \cite{high-toth} said that the number $\frac25$ in
\eqref{toth-ineq} can be replaced by $\frac{2\gamma-1}{1-\gamma}$ and equality
holds only when $n=1$. This means that
\begin{equation}\label{ding}
\frac{1}{2n+\frac1{1-\gamma}-2}\le H(n)-\ln n-\gamma<\frac{1}{2n+\frac13},\quad n\in\mathbb{N}.
\end{equation}
This double inequality was recovered and sharpened in \cite{harm-front, harseq} and \cite[Theorem~2]{qi-cui-jmaa}.
\par
In \cite{yanghar}, basing on an improved Euler-Maclaurin summation formula, some general inequalities for the $n$-th harmonic number $H(n)$ are established, including recovery of the inequality \eqref{toth-ineq}.
\par
In \cite{wuk}, the problems above-mentioned was solved once again by employing
\begin{equation}\label{approx}
H(n)=\ln n+\gamma+\frac1{2n} -\frac12\sum_{i=1}^{q-1}\frac{B_{2i}}{in^{2i}}
-\int_n^\infty \frac{B_{2q}(x)}{x^{2q}}\td x
\end{equation}
and
\begin{equation}
\label{berhar} \int_n^\infty \frac{B_{2q-1}(x)}{x^{2q}}\td
x<\frac{(-1)^qB_{2q}}{2qn^{2q}},
\end{equation}
where $n$ and $q$ are positive integers, $B_i(x)$ are Bernoulli polynomials and $B_{2i}=B_{2i}(0)$ denote Bernoulli numbers for $i\in\mathbb{N}$. For definitions of $B_i(x)$ and $B_{2i}$, please refer to \cite[p.~804]{abram}.
\par
In \cite{wei}, the inequality \eqref{toth-ineq} was verified again by calculus.
\par
In \cite{he-xi}, by utilizing Euler-Maclaurin summation formula, the following general result was obtained:
\begin{equation}
H(n)=\ln n+\gamma+\frac1{2n}-\frac1{12n^2}+\frac12
\sum_{i=3}^{m}\frac{B_{2(i-1)}}{(i-1)n^{2(i-1)}}
+O\biggl(\frac1{n^{2m}}\biggr).
\end{equation}
See also \cite[p.~77]{kuang-3rd}. In fact, this is equivalent to the formula in \cite[p.~259, 6.3.18]{abram}.
\par
In \cite{tothjmaa03}, by considering the decreasing monotonicity of the sequence
\begin{equation}
x_n=\frac1{\bigl|{\sum_{k=n+1}^\infty(-1)^{k-1}\frac1k}\bigr|}-2n,
\end{equation}
it was shown that the best constants $a$ and $b$ such that
\begin{equation}
\frac1{2n+a}\le\Biggl|{\sum_{k=n+1}^\infty(-1)^{k-1}\frac1k}\Biggr|<\frac1{2n+b}
\end{equation}
for $n\ge1$ are $a=\frac1{1-\ln 2}-2$ and $b=1$.
\par
In \cite[Theorem~2.8]{batir-jmaa-06-05-065} and \cite{property-psi-ii.tex}, alternative sharp bounds for $H(n)$ were presented: For $n\in\mathbb{N}$,
\begin{equation}\label{aler-harmonic-ineq}
1+\ln\bigl(\sqrt{e}\,-1\bigr)-\ln\bigl(e^{1/(n+1)}-1\bigr)\le H(n)<\gamma-\ln\bigl(e^{1/(n+1)}-1\bigr).
\end{equation}
The constants $1+\ln\bigl(\sqrt{e}\,-1\bigr)$ and $\gamma$ in \eqref{aler-harmonic-ineq} are the best possible. This improves the result in \cite[pp.~386\nobreakdash--387]{alzer-expo-math-2006}.
\par
In \cite{Infinite-family-Digamma.tex}, it was established that
\begin{equation}\label{new-bounds-harmonic-No}
\ln\biggl(n+\frac12\biggr)+\gamma<H(n)\le\ln\bigl(n+e^{1-\gamma}-1\bigr)+\gamma, \quad n\in\mathbb{N}.
\end{equation}
\par
In \cite{Chao-Ping-AML-10}, it was obtained that
\begin{equation}\label{harmonic-chen-AML-ineq}
\frac1{24\bigl\{n+1/{2\sqrt{6[1-\gamma-\ln(3/2)]}\,}\bigr\}^2}\le H(n)-\ln\biggl(n+\frac12\biggr)-\gamma<\frac1{24(n+1/2)^2}
\end{equation}
for $n\in\mathbb{N}$, where the constants
$$
\frac1{2\sqrt{6[1-\gamma-\ln(3/2)]}\,}
$$
and $\frac12$ are the best possible.
\par
For more information on estimates of harmonic numbers $H(n)$, please refer to \cite{havil, HarmonicNumber},
\cite[pp.~68--86]{kuang-2nd}, \cite[pp.~75--79]{kuang-3rd} and closely-related references therein.
\par
The aim of this paper is to establish a double inequality for bounding harmonic numbers, which is sharp and refines those inequalities above-mentioned.

\begin{thm}\label{harmonic-ineq-ref}
For $n\in\mathbb{N}$, the double inequality
\begin{equation}\label{toth-n^2}
-\frac{1}{12n^2+{2(7-12\gamma)}/({2\gamma-1})}\le H(n)-\ln
n-\frac1{2n}-\gamma<-\frac{1}{12n^2+6/5}
\end{equation}
is valid, with equality in the left-hand side of \eqref{toth-n^2} only when $n=1$, where the scalars $\frac{2(7-12\gamma)}{2\gamma-1}$ and
$\frac65$ in \eqref{toth-n^2} are the best possible.
\end{thm}

\begin{rem}
When $n\ge2$, the double inequality~\eqref{toth-n^2} refines \eqref{harmonic-chen-AML-ineq} and those mentioned before it.
\end{rem}

\section{Proof of Theorem~\ref{harmonic-ineq-ref}}

We now are in a position to prove Theorem~\ref{harmonic-ineq-ref}.
\par
Let
\begin{equation}\label{mfop}
f(x)=\frac1{\ln x+1/{2x}-\psi(x+1)}-12x^2
\end{equation}
for $x\in(0,\infty)$. An easy computation gives
$$
f'(x)=\frac{4x^2\psi'(x+1)-4x+2}{[2x\ln x-2x\psi(x+1)+1]^2}-24x =\frac{4x^2g(x)} {[2x\ln x-2x\psi(x+1)+1]^2},
$$
where
\begin{equation}
g(x)=\psi'(x+1)-\frac1x+\frac1{2x^2}-24x\biggl[\psi(x+1)-\ln x-\frac1{2x}\biggr]^2.
\end{equation}
\par
In~\cite[Theorem~8]{psi-alzer}, the functions
\begin{equation}\label{fn(x)}
F_n(x)=\ln\Gamma(x)-\biggl(x-\frac12\biggr)\ln x+x-\frac12\ln(2\pi) -\sum_{j=1}^{2n}\frac{B_{2j}}{2j(2j-1)x^{2j-1}}
\end{equation}
and
\begin{equation}\label{gn(x)}
G_n(x)=-\ln\Gamma(x)+\biggl(x-\frac12\biggr)\ln x-x+\frac12\ln(2\pi) +\sum_{j=1}^{2n+1}\frac{B_{2j}}{2j(2j-1)x^{2j-1}}
\end{equation}
for $n\ge0$ were proved to be completely monotonic on $(0,\infty)$. This generalizes \cite[Theorem~1]{merkle-jmaa-96} which states that the functions $F_n(x)$ and $G_n(x)$ are convex on $(0,\infty)$. The complete monotonicity of $F_n(x)$ and $G_n(x)$ was proved in~\cite[Theorem~2]{Koumandos-jmaa-06} once again. In particular, the functions
\begin{equation}\label{f1(x)}
\begin{split}
F_2(x)&=\ln\Gamma(x)-\biggl(x-\frac12\biggr)\ln x+x-\frac12\ln(2\pi) \\
&\quad-\frac1{12x}+\frac1{360x^3} -\frac{1}{1260 x^5}+\frac{1}{1680 x^7}
\end{split}
\end{equation}
and
\begin{equation}\label{g1(x)}
\begin{split}
G_1(x)&=-\ln\Gamma(x)+\biggl(x-\frac12\biggr)\ln x-x+\frac12\ln(2\pi)\\ &\quad+\frac1{12x}-\frac1{360x^3}+\frac1{1260x^5}
\end{split}
\end{equation}
are completely monotonic on $(0,\infty)$. Therefore, we have
\begin{multline}\label{psi(x)-double-ineq}
\ln x-\frac{1260 x^5+210 x^4-21 x^2+10}{2520 x^6}<\psi(x)\\
<\ln x-\frac{2520 x^7+420 x^6-42 x^4+20 x^2-21}{5040 x^8}
\end{multline}
and
\begin{multline}
\frac{210 x^8+105 x^7+35 x^6-7 x^4+5 x^2-7}{210 x^9}<\psi'(x)\\
<\frac{210 x^6+105 x^5+35 x^4-7 x^2+5}{210 x^7}
\end{multline}
on $(0,\infty)$. From this, it follows that
\begin{multline}\label{ln-x-psi-upper}
\ln x+\frac1{2x}-\psi(x+1)
=\ln x-\frac1{2x}-\psi(x)\\
<\frac{1260 x^5+210 x^4-21 x^2+10}{2520 x^6}-\frac1{2x}
=\frac{10-21 x^2+210 x^4}{2520 x^{6}}
\end{multline}
and
\begin{align*}
g(x)&>\psi'(x)-\frac1{x^2}-\frac1x+\frac1{2x^2}-\frac{\bigl(10-21 x^2+210 x^4\bigr)^2}{264600x^{11}}\\
&>\frac{210 x^8+105 x^7+35 x^6-7 x^4+5 x^2-7}{210 x^9}-\frac1{x^2}-\frac1x\\
&\quad+\frac1{2x^2}-\frac{\bigl(10-21 x^2+210 x^4\bigr)^2}{264600x^{11}}\\
&=\frac{1659 x^4-8400 x^2-100}{264600 x^{11}}\\
&=\frac{1659 (x-3)^4+19908 (x-3)^3+81186 (x-3)^2+128772 (x-3)+58679}{264600 x^{11}}.
\end{align*}
Hence, the function $g(x)$ is positive on $[3,\infty)$. So the derivative $f'(x)>0$ on $[3,\infty)$, that is, the function $f(x)$ is strictly increasing on $[3,\infty)$.
\par
It is easy to obtain
\begin{gather*}
f(1)=\frac{2 (7-12 \gamma)}{2 \gamma -1}=0.9507\dotsm,\\
f(2)=\frac{4 (48 \gamma +48 \ln2-61)}{5-4 \gamma -4 \ln2}=1.1090\dotsm,\\
f(3)=\frac{3 (108 \gamma +108 \ln3-181)}{5-3 \gamma -3 \ln3}=1.1549\dotsm.
\end{gather*}
This means that the sequence $f(n)$ for $n\in\mathbb{N}$ is strictly increasing.
\par
Employing the inequality~\eqref{ln-x-psi-upper} yields
$$
f(x)>\frac{2520 x^{6}}{10-21 x^2+210 x^4}-12x^2=\frac{12 x^2 \bigl(21 x^2-10\bigr)}{10-21 x^2+210 x^4} \to\frac65
$$
as $x\to\infty$. Utilizing the right-hand side inequality in \eqref{psi(x)-double-ineq} leads to
\begin{align*}
f(x)&=\frac1{\ln x-1/{2x}-\psi(x)}-12x^2\\
&<\frac1{{(2520 x^7+420 x^6-42 x^4+20 x^2-21)}/{5040 x^8}-1/{2x}}-12x^2\\
&=\frac{12 x^2 \bigl(42 x^4-20 x^2+21\bigr)}{420 x^6-42 x^4+20 x^2-21}\\
&\to\frac65
\end{align*}
as $x\to\infty$. As a result, it follows that $\lim_{x\to\infty}f(x)=\frac65$. Therefore, it is derived that $f(1)\le f(n)<\frac65$ for $n\in\mathbb{N}$, equivalently,
$$
\frac{2 (7-12 \gamma)}{2 \gamma -1}\le\frac1{\ln n+1/{2n}-\psi(n+1)}-12n^2<\frac65
$$
which can be rearranged as
$$
\frac{1}{12n^2+{2(7-12\gamma)}/({2\gamma-1})}\ge{\ln n+\frac1{2n}-\psi(n+1)}>\frac{1}{12n^2+6/5}.
$$
Combining this with \eqref{hn-anal} yields \eqref{toth-n^2}. The proof of Theorem~\ref{harmonic-ineq-ref} is proved.


\begin{thebibliography}{99}

\bibitem{abram}
M. Abramowitz and I. A. Stegun (Eds), \textit{Handbook of Mathematical
Functions with Formulas, Graphs, and Mathematical Tables}, 4th printing, with
corrections, Applied Mathematics Series \textbf{55}, National Bureau of
Standards, Washington, 1965.

\bibitem{psi-alzer}
H. Alzer, \textit{On some inequalities for the gamma and psi functions}, Math. Comp. \textbf{66} (1997), no.~217, 373\nobreakdash--389.

\bibitem{alzer-expo-math-2006}
H. Alzer, \textit{Sharp inequalities for the harmonic numbers}, Expo. Math. \textbf{24} (2006), no.~4, 385\nobreakdash--388.

\bibitem{batir-jmaa-06-05-065}
N. Batir, \textit{On some properties of digamma and polygamma functions}, J. Math. Anal. Appl. \textbf{328} (2007), no.~1, 452\nobreakdash--465.

\bibitem{Chao-Ping-AML-10}
Ch.-P. Chen, \textit{Inequalities for the Euler-Mascheroni constant}, Appl. Math. Lett. \textbf{23} (2010), 161\nobreakdash--164; Available online at \url{http://dx.doi.org/10.1016/j.aml.2009.09.005}.

\bibitem{harm-front}
Ch.-P. Chen and F. Qi, \textit{The best bounds of harmonic sequence},
Available online at \url{http://arxiv.org/abs/math/0306233}.

\bibitem{harseq}
Ch.-P. Chen and F. Qi, \textit{The best lower and upper bounds of harmonic
sequence}, RGMIA Res. Rep. Coll. \textbf{6} (2003), no.~2, Art.~14; Available
online at \url{http://www.staff.vu.edu.au/rgmia/v6n2.asp}.

\bibitem{DeTemple}
D. W. DeTemple, \textit{The non-integer property of sums of
reciprocals of consecutive integers}, Math. Gaz. \textbf{75} (1991), 193--194.

\bibitem{havil}
J. Havil, \textit{Gamma: Exploring Euler's Constant}, Princeton, Princeton
University Press, 2003.

\bibitem{high-toth}
R. High, \textit{Asymptotics of the harmonic sum}, Amer. Math. Monthly
\textbf{99} (1992), no.~7, 684--685.

\bibitem{klam}
G. Klambauer, \textit{Problems and Propositions in Analysis}, Marcel Dekker,
New York and Basel, 1979.

\bibitem{Koumandos-jmaa-06}
S. Koumandos, \textit{Remarks on some completely monotonic functions}, J. Math. Anal. Appl. \textbf{324} (2006), no.~2, 1458\nobreakdash--1461.

\bibitem{he-xi}
J.-Ch. Kuang, \textit{Asymptotic estimations of finite sums}, H\=ex\=i
Xu\'eyu\`an Xu\'eb\`ao (Journal of Hexi University) \textbf{18} (2002), no.~2,
1--8. (Chinese)

\bibitem{kuang-2nd}
J.-Ch. Kuang, \textit{Ch\'angy\`ong B\`ud\v{e}ngsh\`i} (\textit{Applied
Inequalities}), 2nd ed., Hunan Education Press, Changsha, China, May 1993.
(Chinese)

\bibitem{kuang-3rd}
J.-Ch. Kuang, \textit{Ch\'angy\`ong B\`ud\v{e}ngsh\`i} (\textit{Applied
Inequalities}), 3rd ed., Sh\=and\=ong K\=exu\'e J\`ish\`u Ch\=ub\v{a}n Sh\`e
(Shandong Science and Technology Press), Ji'nan City, Shandong Province, China,
2004. (Chinese)

\bibitem{merkle-jmaa-96}
M. Merkle, \textit{Logarithmic convexity and inequalities for the gamma function}, J. Math. Anal. Appl.  \textbf{203} (1996), no~2, 369\nobreakdash--380.

\bibitem{PS}
G. P\'{o}lya and G. Szeg\"{o}, \textit{Problems and Theorems in Analysis},
Vol.~I and I\!I, Springer-Verlag, Berlin, Heidelberg, 1972.

\bibitem{qi-cui-jmaa}
F. Qi, R.-Q. Cui, Ch.-P. Chen, and B.-N. Guo, \textit{Some completely
monotonic functions involving polygamma functions and an application}, J.
Math. Anal. Appl. \textbf{310} (2005), no.~1, 303\nobreakdash--308.

\bibitem{property-psi-ii.tex}
F. Qi and B.-N. Guo, \textit{A short proof of monotonicity of a function involving the psi and exponential functions}, Available online at \url{http://arxiv.org/abs/0902.2519}.

\bibitem{Infinite-family-Digamma.tex}
F. Qi, B.-N. Guo, \textit{Sharp inequalities for the psi function and harmonic numbers}, Available online at \url{http://arxiv.org/abs/0902.2524}.

\bibitem{toth-3432}
L. T\'oth, \textit{Problem E 3432}, Amer. Math. Monthly \textbf{98} (1991), no.~3, 264.

\bibitem{tothjmaa03}
L. T\'oth and J. Bukor, \textit{On the alternating series $1-\frac12+\frac13-\frac14+\dotsm$}, J. Math. Anal. Appl. \textbf{282} (2003), no.~1, 21\nobreakdash--25.

\bibitem{wei}
Sh.-R. Wei and B.-Ch. Yang, \textit{A refinement on the Franel inequality},
Zh\=ongy\=ang M\'inz\'u D\`axu\'e Xu\'eb\`ao Z\`ir\v{a}n K\=exu\'e B\v{a}n (J.
Central Univ. Nationalities Natur. Sci. Ed.) \textbf{8} (1999), no.~1, 66--68.
(Chinese)

\bibitem{HarmonicNumber}
E. W. Weisstein, \textit{Harmonic Number}, From MathWorld--A Wolfram Web
Resource. \url{http://mathworld.wolfram.com/HarmonicNumber.asp}.

\bibitem{wuk}
K. Wu and B.-Ch. Yang, \textit{Some refinements of Franel's inequality},
Hu\'an\'an Sh\=\i f\`an D\`axu\'e Xu\'eb\`ao Z\`ir\v{a}n K\=exu\'e B\v{a}n (J.
South China Normal Univ. Natur. Sci. Ed.) (1997), no.~3, 5--8. (Chinese)

\bibitem{yanghar}
B.-Ch. Yang and G.-Q. Wang, \textit{Some inequalities on harmonic series},
Sh\`uxu\'e Y\'anj\=\i u (J. Math. Study) \textbf{29} (1996), no.~3, 90--97.
(Chinese)

\bibitem{young}
R. M. Young, \textit{Euler's constant}, Math. Gaz. \textbf{75} (1991), 187\nobreakdash--190.

\end{thebibliography}
\end{document}